\documentclass[11pt,bibay3]{asl}
\usepackage[legalpaper,margin=1 in]{geometry}
\usepackage[T1]{fontenc} 
\usepackage{hyperref,xspace}
\usepackage[shortlabels]{enumitem}
\newcommand\sectionpadding[1]{\texorpdfstring{}{\hskip #1 \relax}}

\newcommand\Fraisse{Fra\"\i ss\'e\xspace}
\newcommand\HH{\mathcal H}
\newcommand\JJ{\mathcal J}
\newcommand\Nn{\mathbb N}
\newcommand\Qq{\mathbb Q}
\newcommand\Ss{\mathbb S}
\newcommand\Zz{\mathbb Z}
\newcommand\bb{\frak b}
\DeclareMathOperator\Aut{Aut}
\DeclareMathOperator\GL{GL}
\DeclareMathOperator\Sym{Sym}
\newcommand\op{{\textrm{op}}}
\usepackage{wasysym} 
\newcommand\UseBibname[1]{\renewcommand\bibname{\textsc{\normalsize #1}}\renewcommand\refname{\textsc {#1}}\markboth{\bibname}{\bibname}}
\newcommand\etal{\hbox{et al.}} 
\begin{document}
\title{Homogeneity and related topics: \\An Extended Bibliography}

\author{Gregory Cherlin}
\revauthor{Cherlin, Gregory}
\address{Dept.~of Mathematics\\Hill Center, Busch Campus\\Rutgers University\\Piscataway, NJ 08854}
\email{cherlin.math@gmail.com}

\message{Make title ...}
\maketitle

\begin{center}\sc abstract
\end{center}

We give a bibliography of works relating to homogeneous structures in the sense of Fra\"\i ss\'e, and related topics, mainly through 2016, with some narrow updating through 2021.

We first give a list arranged by topics, with references to the main bibliography, which follows.

\medskip
{\it Technical note, 12/2021.} The style file used for the bibliography mishandles long author strings and leaves many 
orphan lines after authors' names.
So the tex file for the bibliography has been edited by hand.
If editing further note the difference 
between \verb|\bibfitem| and \verb|\bibritem| for the first and subsequent occurrences of a group of authors.

\tableofcontents
\newcommand\Subsection[1]{\subsubsection{#1}\indent\par}
\section{\texorpdfstring{\sectionpadding{1.5 em}}{} Literature: a rapid survey}
\label{App:Literature}

The present Appendix contains an outline of some 
literature related in various ways
to the theory of homogeneous structures
and some natural generalizations, mainly through 2916, with a few additions.
Some of this bears
on the classification of
countable homogeneous structures of various types or the close study of their automorphism groups and associated combinatorial problems.

A more focused bibliography accompanies the comprehensive
survey in \cite{Mcp-SHS}.

\medskip

The topics have been organized as follows.

\begin{enumerate}[(I)]
\item  Homogeneous structures:
\begin{enumerate}[(A)]
\item Countable 
\item Uncountable
\end{enumerate}
\item Their Automorphism groups:
\begin{enumerate}[(A)]
\item Algebraic properties
\item Dynamical properties
\item Reconstruction
\item Endomorphisms
\end{enumerate}


\item Generalizations of homogeneity
\begin{enumerate}[(A)]
\item $\aleph_0$-categoricity
\item Transitivity conditions
\item Metric geometry
\item Homomorphism homogeneity
\item Continuous or projective \Fraisse theory
\item Hrushovski amalgamation
\item Random Structures
\item Universality
\end{enumerate}
\end{enumerate}

The boundaries of the subject are very fluid.
Any countable structure can in principle be constructed
by \Fraisse's method---but some  structures 
are most naturally constructed that way, 
and for a number of others this approach is a useful 
complement to other points of view.

\newpage
\begin{enumerate}[(I)]
\item  Homogeneous structures:
\begin{enumerate}[(A)]

\item Countable

\begin{enumerate}[(1)]
\item 
{\it General theory:}\\
\ 
\cite{BaZ-AC, Cam-OPG, Cam-ARS, Cam-IPG, Che-CP,
Fra-RGQ, Fra-PCO, Fra-Ama, Fra-TR1, Fra-TR2,
Hod-MT, Hod-SMT, Jon-URS, Jon-HURS, Mcp-SHS,
Sau-TRA, Tar-MRG}, \cite{TeZ-CMT, Wie-PG, Wie-PGIR}, 
cf.~\cite{Bon-CCnEC, Cam-Tr, Che-CP, CsHMM-CFL,
CuP-AR, DePSS-IRS, DrG-UD, DrG-UDA, Eng-GGT, RoW-US,
Wie-ATPG, Wie-PR};\\
{\it amalgamation bases:} \\
\cite{AlB-BOAB}, 
\cite{Bac-SA, Bac-AIET, BaH-SCA,
Bel-ACAS, Ber-ACDV, Ber-ACMLV, Ber-ABCG, BlG-ROA,
BlG-ABLG, BuFM-FAS, BuFN-APOM, Che-ABCR, Ekl-ACSA,
Fle-AFMAB, GlSW-NOG, Hal-AIS, Hal-ISVAP, Hal-REAS,
Hal-AIRS, Hal-GISA, Hal-FISA, Hal-AGIS, Hal-RSA, HaP-ABFS,
HaS-FBAB, How-ETAS, How-EAS, ImH-AGI*, Lei-ALFp,
Mag-ANG2, Mai-ANG2, Mai-ALFG, Nas-APOM, Ren-EAMS,
Ren-EAR, Ren-FAM, Ren-PAB, Ren-SAM, Sar-ABn2,
Say-ABCA, Say-AUAL, Sho-SAB, Sho-CSAB1, Sho-ABS,
Sho-BSAB, Sho-ABR, Sho-NBAB, Sho-REABR, Sho-RECNB,
Sho-CNBAB, Sho-AFRS, Sho-FISAB, Sho-DREFS, Sho-RSAB,
Sho-FRSAB, Sho-RSABF}.

\item {\it Constructions:} \\
{\it graphs:} 
\cite{ErR-AG, Rad-UG, Hen-HG, Cam-OGHG},\\
cf.~\cite{Nes-AGA}; \\
{\it posets, semi-lattices:} \cite{AlB-ECPO};\\
{\it directed graphs:}
\cite{Hen-2A0, PaPPS-OPO};\\
{\it unbalanced digraphs:}
\cite{EmE-PUD};\\
{\it $n$-fold linear orders:} \cite{Bra-HPS};\\
{\it ultrametric spaces:}
\cite{Bog-MHS};\\
{\it generalized ultrametric spaces:}
\cite{Bra-HPS};\\
{\it metrically homogeneous graphs:}
\cite{Che-2P, DaM-JDE, Mos-EUG, Mos-UGD, Mos-DG}, 
cf.~\cite{Bon-MUG, How-GUIE};\\
{\it Urysohn space:} \cite{Ury25, Ury27a, Ury27b}, cf.~\cite{HuN-FPUS, Kat-UMS};\\
{\it generalized polygons:}
\cite{Ten-FP}, cf.~\cite{Ten-STMT};\\
{\it metric spaces:} 
\cite{Ury25, Ury27a, Ury27b, CaT-LC},
\cite{Hus-UUS, LePRSU-WUS};\\
{\it linear, semi-linear spaces, and Steiner systems:}
\cite{Cam-ILS, Dev-HUS, Dev-USS, Dev-FUSS, DeD-HULS,
Tho-3TSTS, Tre-ISTS};\\
{\it nilpotent groups and rings:} \cite{Bau-FLN, ChSW-HGR};
\\
{\it groups and lattices:} 
\cite{AbT-CHL, Hal-CLFG, Hic-HUG, MaS-ULFG,
Tho-CLF, Tho-CLFL};\\
{\it algebras:} \cite{Gol-TCA};\\
{\it domains, event structures, causal sets:} 
\cite{BoCD-UHSD, Dro-UHES, Dro-FAUD, Dro-UCS, DrG-UDDS,
DrG-UIS}\\
{\it Chu spaces:} \cite{DrZ-BCS};\\
{\it semigroups:} \cite{Ash-ETAB};\\
{\it finite number of countable models:}
\cite{Tan-TC3M, Woo-FCM1, Woo-FCM2};\\
{\it cofinitary permutation groups:} \cite{Ade-TFRH, Cam-PRF, Cam-CPG, Cam-MTSG, Cam-CPGb, Cam-ACPG};\\
{\it probabilistic} \cite{AcFP-IM, AcFNP-IM, AlS-PM,
DoM-PGCL, ErS-PMC}, 
cf.~\cite{ErKR-EKnF, PeV-IM, PrSS-KnF,
PrT-RGPO}.

\item 
{\it Classification:} \\
{\it finite or stable:}
\cite{ChL-SFHS, KnL-SSC, Lac-SHS, Lac-BHS1a, Lac-BH1b, Lac-SH3, Lac-FHD,
Lac-BH2, Lac-HS, Lac-SFH,
LaS-SHB}, 
cf.~\cite{Che-SHS}; \\
{\it finite binary primitive:}
\cite{Che-PBPG, DaGS-BPS, GiS-BPAC, GiHS-BPR1, GiLS-BPL, Wis-PBPG};\\
{\it graphs:} \cite{Gar-HG, LaW-HG, She-SES, She-HG, Woo-4UGT}, 
cf.~\cite{Eno-CHG, Gar-HGS, Gar-HCG};\\
{\it partitioned graphs:} \cite{Ros-H2G}; \index{graph!partitioned}
{\it multipartite edge-colored graphs:} \cite{JeST-HMG, LoT-HCMG};
{\it 3-edge colored complete graph:} \cite{Ara-HS3G, Che-IH3, Lac-SH3, Tri-FH3G};\\
{\it tournaments:} \cite{Lac-HT, Che-HTR};\\
{\it coloured partial orders:} \cite{Sch-HPO, ToT-HCPO};
\\
{\it tournaments with a vertex coloring} \cite{Che-HDG};\\
{\it directed graphs:}
\cite{Che-HDG1, Che-HDG2, Che-HDG, Lac-FHD};\\
{\it bipartite digraphs with partition:} \cite{Ham-2PD};\\
{\it hypergraphs:} \cite{LaT-FH3G},\\
{\it partial orders with vertex coloring} \cite{ToT-HCPO};\\
{\it permutation structures:} \cite{Cam-HP};\\
{\it linear extensions of partial orders:} \cite{DoM-HLPO};\\
{\it ordered graphs \cite{Che-HOGMHG};}\\
{\it finite or locally finite metrically homogeneous graphs:}
 \cite{Cam-6TG, HeP-LFHG, Mcp-DTF} (using \cite{Dun-CUG}), 
 cf.~\cite{Gar-ATG1, Gar-ATG2, Gar-ATG3, Iva-DDRG};\\
{\it infinite metrically homogeneous graphs:}
\cite{AmCM-MH3, Che-2P, Che-HOGMHG};\\
{\it finite homogeneous 3-hypergraphs:} 
\cite{AkL-H3H, LaT-FH3G, Tri-FH3G};\\
{\it homogeneous 3-hypergraphs with one constraint:}
\cite{AkL-H3H}:\\
{\it rings:} \cite{Ber-EQRp, BeC-QERp, BeC-QERpn, BoMP-QESS,
Sar-ECR2, SaW-QEN, SaW-QEp2, SaW-FHRO, SaW-FQE4,
SaW-HFR2n};\\
{\it finite or solvable groups:} \cite{ChF-QEG, ChF-HSG, ChF-HFG};\\
{\it unary algebras:} \cite{Wea-HUDA}.
\item 
{\it Connections with computer science:}\label{Item:ConnectionsCS}\\
{\it cores:}
\cite{Bau-CLP, Bau-CCDG, Bau-FEID, Bod-CCC1, Bod-CCC2, BoP-GBG};
\\
{\it oligomorphic clones:} \cite{BoC-OC, BoJ-CCEI, PeP-PHRS};
\\
{\it constraint satisfaction:}
\cite{Bod-CCS, BoCKO-MCL, BoHM-UACS, BoKM-CSH, BoN-CSCH1, BoN-CSCH2, BoK-CTCS1, BoK-CTCS2,
BoPT-DD1, BoPT-DD2}
\item {\it Model theoretic properties:}\\
{\it generalized metric spaces:} \cite{Con-DGMS, Con-NHMS};\\
{\it homogeneous 3-edge colored complete graphs with
simple theory:} 
\cite{Ara-HS3G};\\
{\it simplicity and supersimplicity:} \cite{AhK-SHS, Ara-OCST,
DeK-G1B, Kop-BP1B, Kop-BSH, Kru-PCC, Pal-GAH};\\
{\it independence relations:} \cite{Con-SIR};\\
{\it no binary homogeneous pseudo-plane:}
\cite{Tho-BHP};\\
{\it finite axiomatizability:} \cite{Lip-FACC, Lip-CCMH, Mcp-FAAC};\\
{\it definable groups:} \cite{Mcp-IG};\\
{\it strongly determined types:} \cite{Iva-SDGC}.
\item {\it Finite approximation and 0-1 laws:}
\cite{BlH-AAG}.
\item {\it Homogenizability:}\\
{\it relational complexity of finite structures:}
\cite{Che-PBPG, ChMS-APG, HaHN-CRS}, 
cf.~\cite{KaK-MSS, KaK-IPG, Sar-WP12, Sar-WP3, Wis-PBPG,
XuGLP-IRAC};\\
{\it relational complexity of infinite structures:}
\cite{Cov-UNFG, Cov-HRS, CoT-HFP, HaHN-BRC}, 
cf.~\cite{Neu-CSTP}.
\item {Applications of Ramsey theory:}\\
{\it reducts:}
\cite{Ben-RBH, JuZ-RQ0, BoP-RRS, BoP-MFRG, BoCP-EPP,
BoPP-RROG, LiP-RGP, Lu-RCCG, PaPPPS-RRPO, Tho-RRG,
Tho-RRH}, 
cf.~\cite{CaLPTW-OARG, Hun-TO, JuZ-RQ0};\\
{\it decidability of positive primitive definability:}
\cite{BoPT-DD1};\\
{\it analysis:} \cite{ArT-RMA}.
\item {\it ``Going forth''} \cite{McL-SGF, McL-FPBF, Vil-HFGF}.
\end{enumerate}

\item {\it Uncountable:}\\
{\it universal locally finite groups:}
\cite{Hic-CULF};
{\it $n$-cardinal spectra} \cite{Ack-nCS}
\end{enumerate}

\newpage
\item Their Automorphism groups:
\begin{enumerate}[(A)]

\item Algebraic properties
\begin{enumerate}[(1)]
\item
{\it Normal subgroups and quotients:}\\
{\it simplicity:} \cite{EvHKL-SGMS, McT-SAG, Tru-IPG1, Tru-SPG, Tru-4G3B};\\
{\it O'Nan-Scott:} \cite{McP-IONS};\\
{\it symmetric group} \cite{AlCM-AQS, Bae-KGAM, Ono-TSI1, Ono-TSI2, Ono-TSI3}, 
cf.~\cite{Ber-CIC, Ber-ECC, Ber-SU, BoF-NSISG, DrG-PCP,
ShT-QSG};\\
{\it $m$-edge colored random graph:} \cite{CaT-AmRG}:\\
{\it linear or semilinear orders:} 
\cite{BaD-NSDT, BlDG-AGOS, BlG-FPOG, BlGGS-AGM,
DrHM-AHSO, DrKT-HSL, GiT-QOPG, GiT-ROS};\\
{\it partial orders:} \cite{GlMR-AGHPO};\\
{\it cycle-free partial orders:} \cite{DrTW-SGCF}, 
\\cf.~\cite{Tru-BCF};\\
{\it trees:} \cite{MoV-AGG};\\
{\it distributive lattices:} \cite{DrM-AUDL};\\
{\it rational topology:} \cite{Tru-CHRW};\\
{\it Urysohn space:} \cite{TeZ-IUS, TeZ-IBU};\\
{\it linear groups:} \cite{Ros-IGL};\\
{\it free amalgamation classes:} \cite{McT-SAG};\\
{\it homeomorphisms:} \cite{And-SGH};\\
{\it multiply transitive actions:} \cite{Cam-NSMT}.
\item 
{\it Maximal subgroups:} \\
{\it symmetric group}
\cite{Bal-MSSG, Bal-ISSG, BaST-MSSG, BrCPPW-MSSG,
CoMM-MISG, CoM-SSGF, McP-MSS, Ric-MSSG}, \\
cf.~\cite{BeS-CSSG};
 \item {\it Small index property:}\\
{\it  general theory:} \cite{Eva-AGIS, Las-PPI, Tru-IPG2};\\
 {\it symmetric group} \cite{DiNT-SISG, Gau-IISG, ScU-PNZ,
 ShT-SISG, Tho-ISG};\\
{\it the random graph:} \cite{Cam-SSIP, HoHLS-SISC, Hru-EPA};\\
{\it Henson graphs:} \cite{Her-EPI, Sol-HLE};\\
{\it linear orders:} \cite{ChT-SI1T, DrT-SIOPG, GlM-AG2H};\\
{\it trees:} \cite{Mol-AGRT, Tru-CFPO};\\
{\it linear groups:} \cite{Eva-SGL, Eva-SIC};\\
{\it relatively free groups:} \cite{BrE-SIF};\\
{\it $\aleph_0$-categorical structures:}
\cite{Her-PISI, HoHLS-SISC};\\
{\it saturated structures:} \cite{Las-ARS, Las-ASS, MeS-SSSI}, 
cf.~\cite{Las-AFM, Las-SIPA, LaS-USSI}.
\item {\it Cofinality:}\\
\cite{DrG-HCD, DrG-CPG, DrHU-ETUC,
DrT-CAGO, DrT-UCDL, Gou-GAQ, McN-SISG, MiS-CUSG,
ShT-UP, ShT-QCS, ShT-UF, ShT-CSSG, Tho-IDCG, Tho-CIPG,
Tho-GDSG}.
\item {\it Bergman property:}\\
\cite{deC-SBG};\\
{\it symmetric group} \cite{Ber-GSG};\\
{\it automorphisms of linear orders:} \cite{DrH-AGC,Mor-1TLO}.
\item {\it Representing words:}
\cite{DoM-RIP, DrT-RWRG, Lyn-WIP, Myc-RIP}.
\item {\it Free subgroups:}
\cite{GaK-UFS}.
\item {\it Embedding theorems:}\\
\cite{Ade-EIPG, BhM-PRF, BhM-LFRG, Bil-GASH, BiJ-RM,
BoDD-MRG, HaKO-SI, Jal-MRT, JaK-RTCT, Mcp-ACC, Mcp-RST,
McW-PGM, Mek-GEAQ, Mel-SCS, Neu-ARW, Nie-USMG,
Nie-UVAG, Tru-EIPG, Usp-IUS}, 
cf.~\cite{Dou-AGUS, Huh-UUMS, MbP-SIUS,
MeSST-RGRW, Mel-GUS, MeS-DFSG, ShT-ISSG}.
\item {\it Regular actions:}\\
\cite{Cam-HCO, CaJ-BG}, 
\cite{CaV-IGU, Dou-GSU}.
\item {\it Orbit growth or profile:}\\
\cite{ApC-OnT, BuH-ELW, Cam-OUS, Cam-OC, Cam-OE, Cam-EMT, Cam-AA, Cam-OC, Cam-SOPG, Cam-CP, CaS-PUS,
CaT-GUS, EnK-RNF, HaR-POH, Mcp-US, Mcp-GRPG,
Mcp-SIPG, Mcp-PGRG, Mer-OnU, Mer-OGOE, Mil-IAG, MnS-LWT,
Nak-GLWT, Pou-RFR, Pou-OAID, Pou-TR, Rob-TLW, Smi-SGR,
Vat-SPC, Vat-PC, Was-CACCT},
cf.~\cite{Pro-CUS};
\item {\it Cycle structure:} \\
{\it primitive groups (esp.~Jordan groups)} \cite{McP-CT}
{\it graphs:} \cite{LoT-CT, Tru-GUG}, 
cf.~\cite{Tru-AAUG};\\
{\it finitary elements:} \cite{BiM-FOAG, Iva-FPT};\\
{\it Parker vectors:} \cite{GeM-PVIG, GeM-PVOG, GeM-CAT}.
\item {\it Isomorphisms up to language:} \cite{CaT-AmRG, Cou-MHG};
\item {\it First order theories:}\\
\cite{GiGT-UAG, Gla-ALO, GlGHJ-2HC, GlMR-AGHPO,
Kni-AGGP, RuS-AGBA}.
\end{enumerate}

\medskip
\item Dynamical properties
\begin{enumerate}[(1)]
\item {\it General theory:}\\
\cite{Aus-MFE, Ell-TD}.
\item {\it Universal minimal flow:}\\
\cite{Bar-MFUS, Ell-UMS, Ell-UMS, KePT-UMF, KeS-AGRP,
MeNT-MUMF, Usp-MCG}.
\item {\it Extreme amenability (and  relation to Ramsey theory):}\\
\cite{Mit-FPM, Usp-EAH};\\
{\it concentration of measure} \cite{GrM-TAII, Led-CMP, Pes-RMEA};\\
{\it Ramsey theory and dynamics:} 
\cite{FaS-RTLG, KePT-RTTD, Moo-ART, MuP-TDUR, Pes-UCTD,
Pes-TGWT, Pes-AC, Pes-mmS, Pes-DIG, Tsa-Hab}\\
cf.~\cite{Vuk-PRG};\\
{\it oligomorphic groups:} 
\cite{Gla-UMSS, GlW-MAGS, GlW-UMGH, Pes-FAMF};\\
{\it homeomorphism groups:} \cite{GlG-UMhH, GlG-MHA};\\
{\it operator algebras:} \cite{GiP-EAG, GiP-OAET}, 
cf.~\cite{Gla-AmRP};\\
{\it metrizable universal flows:}
\cite{BeMT-MUMF, MeNT-MUMF};\\
{\it $\aleph_0$-categorical linear orders:}
\cite{DoGMR-CCLO};
{\it $L^0$} \cite{FaS-RTLG};\\
{\it generic abelian isometry groups:} \cite{MeT-RAG};\\
{\it precompact expansions:} \cite{Ngu-PCE}.
\item {\it Ramsey theory:} \\
\cite{GrRS-RT1, GrRS-RT2, GrRS-RT2p, Hub-BRD, Nes-RT, Nes-RCH,
NeO-Sp, Ngu-UFREA, Ngu-SRT, Sol-SDR}, 
cf.~\cite{Tod-IRS};\\
{\it Ramsey degree in general:}
\cite{Fou-SRT, Fou-SRD, Fou-AR};\\
{\it monotone classes:} \cite{Nes-RCH};\\
{\it canonical partitions:} \cite{LaSV-CPUS, Lar-RTBH};\\
{\it convex equivalence relations:} \cite{Sok-RQ};\\
{\it homogeneous graphs, generalizations:} 
\cite{AbH-MWI, Bod-NRC, Deu-PTG,
Fol-MCS, Hen-EPP, Nes-RCG}
\\
\cite{NeR-PSG, NeR-RTF,
NeR-TPP, NeR-FCS, NeR-PFR, NeR-SRT, NeR-RTS, NeR-RTH,
NeR-RGP, NeR-2RTH, NeR-RCSS, NeR-PRSS, NeR-RON,
NeR-PCRS, PoS-EPRG, PrV-PTPS, Sau-EPTF, Sau-EHP}, \\
cf.~\cite{RoSZ-RFEG, Sau-RFG, Sol-DRT, Spe-EPPT, Spe-RTRT};\\
{\it bipartite graphs:} \cite{FoPS-RDBG};\\
{\it $n$-colorable graphs:} \cite{Ngu-RTnCG};\\
{\it bowtie-free graphs:} \cite{HuN-BF};
\\
{\it trees:} 
\cite{Deu-GRT, Deu-RTRT, DePV-CRT, Fou-SRPT, Mil-RRT,
Mil-RTT, Mil-PTIS, Sol-ART};\\
{\it local order:} \cite{LaNS-PDLO};\\
{\it directed graphs:} \cite{JaLNW-RHDG};\\
{\it partial orders:} \cite{Fou-RDPO, Fou-CRPO, NeR-CPPL, NeR-RPO, PaTW-GO, Sok-RPFP,
Sok-RPFP2}, 
cf.~\cite{HuN-FPH}; \\
{\it boron trees:} \cite{Jas-HRP, Jas-RBT};\\
{\it metric spaces:} \cite{DiP-MR, ErGMRSS-ERT2,
ErGMRSS-ERT3, Jas-HRP, Nes-RTM, Nes-MSR, Ngu-UUS,
Ngu-SRTD};\\
{\it matroids:} \cite{NePT-AMA};\\
{\it vector space:} \cite{LaNPS-PIVS, GrLR-RTC1, GrLR-RTC2,
Spe-RTS};\\
{\it affine space:} \cite{NePRV-COT};\\
{\it inner product spaces:} \cite{Jas-HRP};
\\
{\it Steiner systems:} \cite{BhNRR-RSS}\\
{\it cubes:} \cite{NePRV-COHJ1, NePRV-COHJ2, Pro-NOCC,
Sol-RTS};
\\
{\it indivisibility:} 
\cite{DeLPS-IUM, ElZS-IKn, ElZS-RRS, ElZS-DHDG,
ElZS-DHDG, ElZS-DHH, ElZS-GVP, ErHP-SEG, KoR-CUG,
LoAN-OSUS, Mel-GDUS, Ngu-BRD, NgS-USOS, NgS-WIMS,
Pou-RI, Sau-VPP, Sau-CVP, Sau-AWI},\\
cf.~\cite{ErH-FIC, ImKT-DIG, ImSTW-2DGG, Sau-EHP,
Sau-FSCS, Sau-CSRG, SaWR-CRT, WaZ-DLFT};
\\
{\it pigeonhole property:} \cite{BoCD-TOPP, BoD-PHP};\\
{\it inexhaustible structures:}
\\
 \cite{BoD-IG, BoSZ-IHS};
{\it affine and projective space:} 
\cite{DeV-PAPR};\\
cf.~\cite{GrRS-RT1, GrRS-RT2, GrRS-RT2p};
cf.~\cite{LaNS-NHRS}.
\item {\it Amenability and unique ergodicity:}\\
{\it equivalence relations:} \cite{Iva-AAG, PaS-UEHDG}.
{\it graphs and digraphs:} \cite{AnKL-UEAG, Zuc-AUE}.
\item {\it EPPA (Hrushovski property), ample generics, generic automorphisms:}\\
\cite{Con-HMS, Con-PIGMS, EvHKN-2G_AM, Her-EPI, HeL-EPA, 
Hod-LGFMP, HoO-FCHC, Iva-GE,
Iva-AHS, Iva-SBAG, Iva-GSDO, KeR-TAGA, KuT-GAPO,
LoT-GEHS, McT-CCC, Ros-GI, Ros-FAGA1, Ros-FAGA2}, 
cf.~\cite{Pes-HSV, PeU-RFG, Slu-NGP, Sol-EPI, Tru-GAHS,
Ver-GPI}.
\item {\it Strong non-local compactness:} \cite{Mal-TOS}.
\item {\it Measures on models:} \cite{Alb-MRG}.
\end{enumerate}
\medskip

\item {\it Reconstruction (see also small index property):}\\
\cite{Bar-AGOC, Bar-RCG, BaM-RHRS, BeR-RFH, GoR-GQA,
KuR-RUS, LeR-RLCS, Ros-ACGH, RoS-ACFP, Rub-ABA,
Rub-RBA, Rub-RTS, Rub-RBAAG, Rub-RTAG, Rub-RCC,
Rub-LMG, Rub-RFM, RuR-EH1O, Slu-ACFP, Tru-QAG,
Tru-AGCO}.
\medskip

\item {\it Endomorphisms:}\\
{\it maximal subgroups:} \cite{McP-EFL};\\
{\it random graph:} \cite{DeD-EMRG, Dol-EMRP};\\
{\it partial orders:} \cite{Mas-EMPO};\\
{\it embedding theorem} \cite{DoM-EMU};\\
{\it Bergman property:} \cite{Dol-BPEM}.
\end{enumerate}

\newpage
\item Generalizations of homogeneity
\begin{enumerate}[(A)]
\item {\it $\aleph_0$-categoricity:} \\
{\it general theory:} \cite{Eng-UFM, RyN-CC, Sve-CC}, 
cf.~\cite{Car-FCTP, Che-ACC, ClK-RHS,
Eva-ECC, Grz-LUCC, Grz-DPCC, Hau-VSC, Sch-CCIS,
Wea-CCT1, WeL-CCT2},
\cite{Cam-OPG, Cam-PGHC, Cam-PGr, Cam-IPG, Cam-OPG2, DiM-PG, Eva-PGMS, Hod-CPG};\\
{\it constructions:} \cite{Ash-UCC},
\cite{CaPZ-GCM, Ehr-CCCT, Gla-2A0},
cf. \cite{Iva-CCC},  
\cite{Tho-CCEP, Per-FNCM, Woo-3IT};\\
{\it total categoricity:} \cite{Ahl-TCMT, Ahl-ASM, AhZ-QFA, AhZ-Z/4Z,
AhZ-IS, BaCM-TCGR, Che-TCS, ChHL-CSS, Hod-STC, Hru-TC,
Sch-CICC, Ten-TCG, Zil-SMCC1, Zil-TCT1, Zil-TCCG, Zil-SMCC2,
Zil-SMCC3},\\
cf.~\cite{HoHM-CCRC, Iva-DLAC, Iva-NGC, IvM-AE, Kos-SFC,
Vas-CCDAG};\\
{\it covers:} 
\cite{ChHS-CLG, Eva-SFC, Eva-CGFC, Eva-FCFK, Eva-FCCC,
EvG-KCG, EvH-CRC, EvH-CCPG, EvH-AGFC, EvIM-FC,
EvP-2CFC, EvR-BFC, HoP-CS, Iva-CP, Iva-FCCH, IvM-SDT,
Pas-AFFC, Pas-FCG, Pas-NHC};\\
{\it graphs:} \cite{Shi-CCG, Whe-CUT, Whe-TNCG};\\
{\it colored linear orders:} \cite{MwT-CLO, MwT-FCLO, Ros-CCLO,
Ros-LO}, 
cf.~\cite{CrT-AS, CrT-QAS, Gla-OPG, HeMMNT-CCWo,
Kul-WCM, Kul-CCoM, KuM-MCO};\\
{\it Boolean algebras with ideals:} 
\cite{Ala-CCBA, Pal-BADI}, 
cf.~\cite{Hei-WFBA};\\
{\it multitrees (``reticles'')} 
\cite{Pun-MCD, Pun-MCMT};\\
{\it partial orders:} \cite{Pou-EOU2C, Sch-CCPO};\\
{\it distributive lattices:} \cite{Sch-CCDL};
\\
{\it rings:} 
\cite{BaR-CCSR, Che-CCN1, Che-CCN2,
MaR-CCRB, Ros-PQER, Ros-GL2B, Ros-CCG};
\\
{\it groups:}
\cite{App-BPG, App-CCG, App-ECCG, App-FEBP,
App-ECCG, ArM-SCCG,
BaCM-TCGR,
ChR-AFG,
Fel-SCCG, Fel-CCSG, Iva-CCG,
Iva-DIT, IvM-AUG,
Mcp-AUCCG,
Ros-PQER, Ros-GL2B, Ros-CCG,
SaW-PENG, SaW-QEe4, Wil-SCCG}
cf.~\cite{IvM-ECCG};\\
{\it quasi-groups:} \cite{Shi-CQQ};\\
{\it bilinear maps:} \cite{Bau-CCBM};\\
{\it quasi-varieties:} \cite{BaL-UHC, Pal-CQV, Pal-DCQV,
Pal-CPHT, Pal-CPHC, Pal-CHC};\\
{\it e.c.~structures for some universal Horn classes:} \cite{Alb-PEC};\\
{\it automorphism groups:} \cite{Bar-AGOC, Tsa-UOG};\\
{\it Keisler measures:} \cite{Ens-AIM};
\\
{\it model companions:} \cite{Sar-MCCC};\\
{\it simplicity:} \cite{Pal-CCST, PaW-SCMT};\\
{\it orbit growth:} \cite{Pal-AFCC}, cf.~\cite{BoT-HMPG};\\
{\it computable models, decidability:} 
\cite{KhM-CCCTA, Mor-CCDM, Puz-CCT, Sch-DCCT,
Sch-CCT, Sch-RNF, Sch-TCCPO, Stu-HFS};\\
cf. \cite{Per-CCFA, Ric-DS}.
\begin{enumerate}[(1)]
\item {\it Coordinatized by indiscernible sets:} \cite{Lac-IS}
\item {\it Tree decomposable:} \\
\cite{Lac-CCI1, Lac-CIG, Lac-CCI2, Lac-TDS}, cf.~\cite{Pal-CCUT}.
\item {\it Smooth approximation:} \\
\cite{Che-LFFT, ChH-FSFT, Hru-FSFT, KaLM-SAF} 
\item {\it Simple theories:} \cite{Ara-OCST, EvW-SCCG}.
\end{enumerate}
\medskip

\item Transitivity conditions
\begin{enumerate}[(1)]
\item {\it Setwise homogeneous:}\\ 
{\it in general:} \cite{Mcp-HIPG}, 
cf.~\cite{LoM-OEPG};
\\
{\it graphs:} \cite{DrGMS-SHG, DrGM-EO, GrMPR-SHDG};
\\
{\it directed graphs:} \cite{GrMPR-SHDG}.
\item {\it $k$-homogeneity and variants:} \\
\cite{Bro-WnH, Cam-LW, Cam-TUS, Cam-PGUS, Cam-OUS, Cam-OUS2, Cam-OUS3, Cam-OUS4, Dro-kHUP, Dro-kHRT, Hig-HR,
Hug-kHG, Kan-4HG, Kan-kHG, Kie-HRBP, Mcp-HIPG, Mar-PGT,
Neu-HIPG, Tru-PHO, Wie-kHPG, Yos-kHPG}, 
cf.~\cite{Haj-HIPG, Ros-NCC, ShT-HIPG};\\
{\it graphs:} \cite{FaLLP-L2ATB, GoK-kHG, Gar-SCG, LiW-TPOS},
cf.~\cite{Dej-K4K222HG, Dej-C4UHG, Dej-OSDT, Gra-kCSG,
IsJP-KUH, LiSS-sAT, Ron-HG, ShS-HPG, YaF-WsATG};\\
{\it linear orders:} \cite{DrS-OAG};\\
{\it circular orders:} \cite{CaT-1TCO, GiH-OhS, KuM-MCO};\\
{\it partial orders:} \cite{Dro-POT1, Dro-POT2, DrM-kHPG,
DrMM-UHPO, SaW-PHPO};\\
{\it cycle-free partial orders:} 
\cite{CrTW-kCST, GrT-CFPO, Tru-kCST, War-kCSCF}, 
cf.~\cite{Mol-EG1, Mol-EG2};\\
{\it linear  spaces:} \cite{Dev-dHLS};\\
{\it affine or projective space:} \cite{Tho-GPS1, Tho-GPS2}, 
cf.~\cite{Tho-ISG};\\
{\it real measurement:} \cite{Alp-RMS, Alp-OHFU};
\item {\it Canonical expansions:}
\cite{Tho-ECCS}
\item {\it $1$-homogeneity in the sense of Myers:}\\
\cite{GaM-SP1H, McA-HCP, Mye-nHG, Mye-1HG}.
\item {\it Distance transitive graphs:}\\
{\it finite:} \cite{Bon-ADTG, Bon-FPDT, Cam-MTG, Cam-6TG, Gar-SG,
Iva-DTG};\\
{\it infinite:} \cite{Cam-Cen, Mcp-DTF, Mol-DTIG}, 
cf.~\cite{Mol-LFG};\\
{\it imprimitive:} \cite{AlH-ST, Smi-PIG};\\
{\it distance transitive with more than one end:}
 \cite{HaP-TCIG}, using \cite{DuK-VC};\\
{\it distance regular:} \cite{BrCN-DRG}.
\item {\it Highly arc-transitive digraphs:}\\
\cite{Ama-DHAT, AmT-HATD, AmT-CFHA, CaPW-HATD,
Che-HATD, DeMS-HATD, MaMSZ-HATD, MaMMSZ-HATD,
Mol-DATD, Neu-HATD, Pra-HATD}, \\
cf.~\cite{DeJLP-LGGT, GiLP-sATG, GiLP-3ATG, GiLP-5ATG,
JiDLP-GTG, Mol-LCG}.
\item {\it Descendant-homogeneous digraphs:}\\
\cite{AmET-CDHD, AmT-DHD}.
\item {\it Homogeneous with respect to connected induced subgraphs or digraphs:}\\
\cite{GrM-CHD, Ham-LFCHD, Ham-ETG, HaH-CHDE, HaP-TCIG},
cf.~\cite{DuK-VC, KrM-ME, KrM-QIGT}.
\item {\it Orbit-homogeneity:}
\cite{CaD-OH}.
\item {\it Jordan groups:}\\
{\it survey:} \cite{BhMMN-IPG, Mcp-SJG}; 
\cite{Ade-GJG, Ade-STSS, Ade-IJPG, Ade-IIJG, AdM-CJG,
AdN-PJS, AdN-IBPG, AdN-RRB, BhM-JGB, Eva-HG, Hic-JRFP,
Hru-UMS, Hyt-LMG, Joh-STS, Kan-HDGL, McD-JG,
Neu-PPG, Neu-IJG}
\item {\it Universal transversal property:}\\
\cite{ArC-GHG}.
\item {\it fine partition:} \cite{HoM-RSFIS, Lac-UE}.
\item {\it Extension properties:} \\
{\it graphs:} \cite{Ana-APGPG, AnC-GPAP, AnC-APPG, AnC-GSAP,
AnC-GPAP1, AnC-GPAP2, AnC-CQPnEC,
BaBB-3ECAP, BaBS-nECAP, BaBMP-nECRD,
BlEH-PGA, BlR-EGEP, Bon-nEC, BoC-APG,
BoC-2ECLC, BoC-2ECGT, BoC-APT, BoC-APGE,
BoHK-SR3EC, CaS-SRnEC, DaM-JDE, ErHK-FSU,
ErP-SIS, Exo-APG, ExH-GCAP, ExH-SGAP,
Fag-PFM, RoSW-kUG}, \\
cf.~\cite{Fla-PCSR, Tro-UMS};\\
{\it triangle free graphs:} 
\cite{AlCH-TFA, EZL-3EC, Lar-TFGS, Pac-ISCN, PaS-2SU}, 
cf.~\cite{AlR-TFGS};
{\it tournaments:} \cite{GrS-CSTP}.
\item {\it Association schemes:} \cite{AlBC-ASPG}.
\end{enumerate}

\medskip
\item {\it Metric geometry:}\\
\cite{Bir-MFG, Bog-MHS, Bus-LDP, Bus-SFS, Bus-MMFS,
Bus-2PG, Bus-LMG, Bus-GG, Bus-AG, BuP-DG,
DaW-MHR, Fre-RHL1, Fre-RHL2, Nag-TGNC, Nag-HSB,
Nag-WTF, Tit-CEM, Tit-CEMn, Tit-EHGL, Tit-TGM,
Wan-2MS, Wan-2PH}, 
cf.~\cite{Tit-G3T, Tit-G2TC}.
\medskip

\item {\it Homomorphism-homogeneity:}\\
\cite{CaL-PHH, CaN-HHS, DoM-HHL, DoJ-HHP, HaHM-HHLG,
IlMR-FHHT, IlMR-HHG, JuM-HHMA, Loc-CHHG, LoT-HH,
Mas-HHPO, Mas-HHTL, Mas-HHFA, Mas-HHPLG, Mas-HHOGL,
MaNS-HHBR, MaNS-FHHB, MaP-HHS, RuS-HHG},\\
cf. \cite{PeP-PHRS}.
\medskip

\item {\it Continuous and projective \Fraisse theory:}\\
{\it continuous \Fraisse constructions:} 
\cite{AvSCGM-BSUD, BaK-LF, BeM-USECS, BeY-FLM,
Cam-QPFL, GaK-GS, IrS-PAMT, IrS-PFPA, Kub-FS, KuS-UGS,
Sch-FTMS, Usv-GSMS}, \\
cf.~\cite{ArB-UEC, Dza-UEC, Dza-UUEC, Kub-IRFL, Mas-RUS,
ShZ-UMKS};
\\
{\it generic automorphisms, continuous case:} 
\cite{BeBM-PTG, GuI-PGRT, Hod-PAQ, KaL-PGAG,
Kwi-HCAG, Kwi-LCPA, RiZ-PTFG};\\
{\it Ramsey theory:} \cite{Kai-ARTM};\\
{\it linear metric rigidity:} \cite{MePV-EMLR, MePV-LRMS};\\
{\it background:} \cite{BeBHU-MTMS};\\
{\it neostability:} \cite{CoT-PUS}.

\medskip
\item {\it Hrushovski amalgamation:}\\
\cite{AnI-SSG, Are-HUG, Bal-ASMP, Bal-LBCI, Bal-RHS, Bal-EG, BaH-CR2F, BaH-CSCR, BaH-CSRk, BaH-CSMC, BaI-GPP, BaS-DFG, Bau-AC, Bau-UCG, BaHMW-BF, BaMZ-FVS, BaMZ-HF, BaMZ-RF, BaP-FPS,
EaO-CATF, Eva-HAD, Eva-SIPG, Eva-CCP,
Eva-BTSS, Eva-TNT, EvF-GHC1, EvF-GHC2, EvGT-SAG,
EvP-GAC, Goo-HG, Has-IFMR, Has-HAC, HaH-FSL, HaH-DMP,
Her-SGS, Her-STFT, Hol-GF, Hol-FVS, Hol-MCSM, Hru-CCP,
Hru-SMF, Hru-NSM, Hru-SLG, Ike-MFD, Ike-GPP, Ike-SSGG,
Ike-AIG, IkK-SGS, IkKT-GSSA, KuL-GS, PiT-ACC,
Poi-CE, Poi-EC, Poi-AH,
Pou-SCRT, Pou-SGS, Pou-SFGS, Pou-SFHMT, PoW-SPRT,
Sud-VGC, Sud-GGIW1, Sud-SGM, Sud-GGIW2, Sud-GL,
Ten-HnG, Tsu-RAST, Tsu-AT,
Ver-EISM, VeY-CMT, Wag-RSD, Wag-HA, Zie-FFMR,
Zie-NSM, Zil-CSAV, Zil-PEACF, Zil-CQEC, Zil-CMG};\\
{\it relation to random structures:}
\cite{Bal-NMC01, Bal-FIMT, Bal-MRG, BaM-DTLL, BaS-RSG,
BaS-SG, Bey-RSPF, Bey-RHPF, BeH-ACPF, DeN-ASnG,
DoL-ORS, Lyn-ASU, Lyn-PUF, Lyn-PSRG1, Lyn-PSRG2, Lyn-E01,
Lyn-PGMT, Lyn-CLRG, Lyn-CLRGD}.
\medskip

\item{\it  Random Structures:}\\
\cite{AlF-RGO, Boh-TFP, BoPLPSSV-DT, BoK-HFP, BoK-DCTF,
Bol-RG, BoJ-IRGG, BoJW-nOG, BrT-RSC, Cam-RG, ChS-LRG,
Com-TTFM, Com-LAAC1, Com-CDP, CoHS-IAP, Com-AP1PO,
Com-LAAC2, Com-LLC, DrG-RpG, DrK-RRS, Erd-RTG, Erd-CCG,
ErR-AG, Fag-PFM, Fag-MTPP, GlKLT-FSPC, Gra-AAFS,
Las-SS01, ShS-01SRG, Spe-NDM, Spe-PM, Spe-10L, Spe-PM2,
Spe-EEA, SpJ-T01, Tar-MRG, Ver-RMSU, Win-ROFD, Win-RO,
Win-RODk, Win-TRO, Win-RS01},
cf.~\cite{KoPR-KnF, KoPR-KnF2};\\
\item Theory of relations (general structures)
\cite{McPW-CSA}.
\medskip

\goodbreak
\item Universality
\begin{enumerate}[(1)]
\item {\it Countable case:} \\
{\it graphs:} 
\cite{ArBM-CAC, BoT-GCSG, Bon-PHGC, Bon-HA,
ChK-UPF, ChSh-UFT, ChSS-UFS, ChS-FC,
ChS-FSFS, ChST-OBT, ChT-FP2B, ChT-2CP,
FuK-CUG, FuK-UGT, Kom-UG, KoMP-UG, KoP-UBS,
KoP-UE, Mos-EUG, Nur-GMFC, Rad-UGUF, Rad-UG}, 
cf.~\cite{BrHM-UGHP, Kom-OSUG, Pac-MPCG};\\
{\it width-2 orders:} \cite{BoD-W2O};\\
{\it partial orders:} \cite{HuN-UGPO};\\
{\it rings:} \cite{Sar-UCR, SaW-UCR};\\
{\it permutation patterns:} 
\cite{AtMR-PAS, HuR-UPC, Bon-CP2, HuR-UPC}, cf.~{MaT-EPM};\\
\item {\it Uncountable case:} \\
{\it structures with $n$-dimensional amalgamation:} \cite{Mek-US1};
\\
{\it graphs:} \cite{KoP-UE, DzS-UGSS, KoP-UE, She-CST, She-UGCH,
She-UGCHR, She-US, She-UGOB, Tho-SUF};\\
{\it groups:} 
\cite{GoSW-ULNG, GrS-ULFG, Hic-CULF, She-USAG, She-NEU,
She-NEUAG, ShS-KPTF, ShS-KPTFR};\\
{\it topological groups:} \cite{Shk-UATG};\\
{\it linear orders:} \cite{KoS-UO, Moo-UAL, She-IR};\\
{\it partial orders:} \cite{Joh-UIPO};\\
{\it cardinal spectra} \cite{Ack-nCS};\\
{\it bipartite graphs:} \cite{GoGK-BU};\\
{\it topological spaces:} 
\cite{MaNO-USRT, MaT-URS, Tod-FS2}; \\
{\it metric space:} \cite{Kat-UMS};\\
{\it Banach spaces:} 
\cite{BrK-UBSC, BrK-UBSE, BrK-IUBS, Kos-UMUF, Kos-UOBS,
ShU-BSGOP, Szl-NESRU};\\
{\it von Neumann algebras:} \cite{Oza-UII1};\\
{\it using club guessing} \cite{Dza-CGU};\\
{\it universal models:} \cite{Dza-UM, DzS-EU, DzS-NEU, KoS-USUS}.
\end{enumerate}
\end{enumerate}
\end{enumerate}

\clearpage
\phantomsection
\newcommand\LongAuthors[1]{\vspace{\intergroupsep}\noindent 
\textsc{#1}}
\UseBibname{\large Bibliography on Homogeneity}
\bibliographystyle{asl}

\end{document}